\documentclass{article}
\begin{document}

\pagestyle {myheadings}

\markright {Marcelo Gomez Morteo}

\title{On the Semisimplicity of the Action of the Frobenius on Etale Cohomology}
\author{Marcelo Gomez Morteo}

\maketitle
\vspace{16pt}

\begin{abstract}

\vspace{16pt}

 We give a proof of the semisimplicity of the action of the geometric frobenius on etale cohomology. The proof is based on [MGM11] and on the Weil Conjectures, ie on the Riemann Hypothesis for non singular projective varieties over finite fields.

\vspace{16pt}

 Algebraic and Topological K Theory l-adic Completion of Spectra.
\end{abstract}

\textbf{Introduction:}

\vspace{16pt}

 We are going to work with the abelian category $\cal B$${(l)}_{*}$ See the definition below and also [B83]. In [B83] Bousfield defined an universal functor $\cal U$$: Z_{(l)}-modules \mapsto $$\cal B$${(l)}_{*}$ which will be crucial here. We start with the isomorphism ([T89])

\vspace{12pt}

(1) \[ \pi_{0}([L_{E(1)}K(X_{\infty})]^{l})\otimes Q  \simeq \oplus_{i=0}^{d}H^{2i}_{et}(X_{\infty},Q_{(l)}(i))\]

Here $L_{(E(1)}$ is is the Bousfield localization of $E(1)$ where $E(1)$ is such that $\cal K$$\sb{(l)}=\vee_{i=1}^{l-2} \Sigma^{2i} E(1)$.(See below the definition of $E(1)$ and also [B83]) $\cal K$$_{(l)}$ is the $l$ localized topological K spectrum $\cal K$. In (1), $K(X_{\infty})$ is the Quillen's algebraic K theory spectrum on $X_{\infty}=X\otimes F$, where $F$ is the algebraic closure of the finite field $F_{q}$ with $q=p^{s}$ and $p \neq l$. $X$ is a nonsingular projective variety over $F_{q}$ with finite dimension and $d$ is the dimension of $X_{\infty}$ . $L_{E(1)}$ denotes $E(1)$ localization, and $[L_{E(1)}K(X_{\infty})]^{l}$ denotes the $l$-adic completion of the spectra $L_{E(1)}K(X_{\infty})$.Here $H^{2i}_{et}(X_{\infty},Q_{l}(i))$ are the etale cohomology groups with coefficients in $Q_{l}$,which are $Q_{l}$ vector spaces.

\vspace{8pt}

We show that $\cal U$$(\pi_{0}([L_{E(1)}K(X_{\infty})]^{l})\otimes Q) \simeq [E_{0}(1)(K(X_{\infty}))]^{l} \otimes Q$ (2)

\vspace{8pt}

where on the right hand side of the equality, $[E_{0}(1)(K(X_{\infty}))]^{l}$ is the $l$-adic completion of the $Z_{(l)}$-module $E_{0}(1)(K(X_{\infty}))$ which is an object in $\cal B$ (See the definition of $\cal B$ below). This isomorphism is dependent on highly non trivial facts, such as the fact that $\pi_{0}L_{E(1)}K(X_{\infty})$ is an $l$-reduced group, which follows from [MGM11]. Being $l$-reduced group means:

\vspace{6pt}

\[Hom(Q/Z_{(l)},\pi_{0}(L_{E(1)}K(X_{\infty})) )=0\]

\vspace{6pt}

We then obtain the isomorphism: $[(E_{0}(1)K(X_{\infty}))]^{l} \otimes Q \simeq $$\cal U$$(\oplus_{i=0}^{2d}H_{et}(X_{\infty},Q_{(l)}(i)))$. Therefore the functor $\cal U$ allows to study the action of the geometric frobenius $\Phi_{X}:X \mapsto X$ on $\oplus_{i=0}^{2d}H_{et}(X_{\infty},Q_{(l)}(i))$ through the action of $\cal U$$(\Phi_{X})$ on $[E_{0}(1)(K(X_{\infty}))]^{l} \otimes Q$. We prove that $(E_{0}(1)(K(X_{\infty}))) \otimes Q$ is dense in $[E_{0}(1)(K(X_{\infty}))]^{l}\otimes Q$ with the $l$-adic topology and finally show that $\cal U$$(\Phi_{X})$ is an Adams operation on $[E_{0}(1)(K(X_{\infty}))]^{l} \otimes Q$  because of the Weil Conjectures and we conclude that the action of $\Phi_{X}$ on the etale cohomological spaces $H^{2i}_{et}(X_{\infty},Q_{(l)}(i))$ is semisimple.

\vspace{16pt}
\emph{ The Category $\cal B$${(l)}$${}\sb{*}:$ }

\vspace{16pt}

We begin by describing an abelian category, denoted $\cal B$${(l)}$${}\sb{*}$,equivalent to the category of $E(1)_{*}E(1)$-comodules (see [B83], 10.3) Bousfield describes $\cal B$${(l)}$${}\sb{*}$ as follows: Let $l$ be an odd prime and let $\cal B$ denote the category of $Z_{(l)}[Z_{(l)}^{*}]$-modules for the group ring $Z_{l}[Z_{l}^{*}]$, where $Z_{(l)}^{*}$ are the units in $Z_{(l)}$,with the action by the group ring defined by Adams operations $\Psi^{k}:M \mapsto M$ which are automorphisms and satisfy the following:

\vspace{16pt}

i) There is an eigenspace decomposition

\vspace{20pt}

\[ M\otimes Q \cong \bigoplus _{j\in Z} W_{j(l-1)}\]

\vspace{20pt}

such that for all $w\in W_{j(l-1)}$ and $k\in Z_{(l)}$,

\vspace{25pt}

\[(\Psi^{k}\otimes id)w=k^{j(l-1)}w \]

\vspace{20pt}

ii) For all $x\in M$ there is a finitely generated submodule $C(x)$ containing $x$, satisfying:
for all $m \geq 1$ there is an $n$ such that the action of $Z_{(l)}^{*}$ on $C(x)/l^{m}C(x)$ factors through the quotient of $(Z/l^{n+1})^{*}$ by a subgroup of order $l-1$.

\vspace{16pt}

To build the category $\cal B$${(l)}$${}\sb{*}$ out of the above category $\cal B$, we additionally need the following:

\vspace{1pt}

Let $T^{j(l-1)}:\cal B$ $ \mapsto \cal B$ with $j\in Z$ denote the following equivalence:

\vspace{1pt}

For all $M$ in $\cal B$, $T^{j(l-1)}(M)=M$ as $Z_{(l)}$-module, but not as $Z_{(l)}[Z_{(l)}^{*}]$-module since the Adams operations in  $T^{j(l-1)}(M)$ are now $k^{j(l-1)}\Psi^{k}:M \mapsto M$ where $\Psi^{k}$ is the Adams operation of multiplication by $k$ in $\cal B$.
Now an object in $\cal B$${(l)}$${}\sb{*}$ is defined as a collection of modules $M=(M_{n})_{n\in Z}$, with $M_{n}$ in $\cal B$ together with a collection of isomorphisms for all $n\in Z$,

\vspace{25pt}

\[T^{l-1}(M_{n}) \mapsto M_{n+2(l-1)} \]

\vspace{20pt}

Note that the category $\cal B$ can be viewed as the subcategory of $\cal B$${(l)}$${}\sb{*}$ consisting of those objects $(M_{n})_{n\in Z}$ such that $M_{n}=M$ if $n$ is congruent to $0$ mod $2(l-1)$ and $0$ otherwise

\vspace{25pt}

In [B83] Bousfield constructs a functor $\cal U:$ $\pi_{*}(E(1)-Mod \mapsto $ $\cal B$ ${(l)}_{*}$. For $H \in \pi_{*}(E(1)-Mod $, let $\cal U$ in $\cal B$ consist of the objects $\cal U$ $(H_{n})$ in $\cal B$ for all $n \in Z$.

\vspace{20pt}

\emph{ The Spectrum $E(1)$ and its homology theory $E(1)_{*}$:}

\vspace{16pt}

Given $E(1)$, which by construction depends on the prime $l$, there is a map $ E(1) \mapsto $$\cal K$$_{l}$ which is a ring morphism (see [R] Chapter VI Theorem 3.28) and verifies the equivalence $\cal K$$\sb{(l)}=\vee_{i=1}^{l-2} \Sigma^{2i} E(1)$. There are Adams operations $\Psi^{k}:E(1) \mapsto E(1)$ with $k$ in $Z_{l}^{*}$ which are the units in $Z_{l}$.These Adams operations are ring spectra equivalences and $\Psi^{k}$ carries $\nu^{j}$ to $k^{j(l-1)}\nu^{j}$ in $\pi_{2j(l-1)}E(1)$ for each integer $j$ where $\nu$ is such that $\pi_{*}E(1)=Z_{(l)}[\nu,\nu^{-1}]$ and $\nu$ has degree $2(l-1)$.
\vspace{1pt}
Another property of $E(1)$ is that $E(1)$ localization is the same as $ cal K$ $_{(l)}$ localization.

\vspace{16pt}

The homology $E(1)_{*}(X)$ with $X$ a spectrum also has Adams operations $\Psi^{k}: E(1)_{*}(X) \mapsto E(1)_{*}(X)$. One checks that $\Psi^{k}(\nu^{j}x)=k^{j(l-1)}\nu^{j}\Psi^{k}(x)$ for each integer $j$ and $k$ in $Z_{(l)}^{*}$ and $x\in E(1)_{*}(X)$. The multiplication by $\nu^{j}$ induces an isomorphism $\nu^{j}:T^{j(l-1)}E(1)_{n}(X) \mapsto E(1)_{n+2j(l-1)}(X)$ in $\cal B$${(l)}$${}\sb{*}$ for each $j,n \in Z$. It follows that $E(1)_{*}(X)$ is in $\cal B$${(l)}$${}\sb{*}$  for each spectrum $X$ in $\cal S$ by taking $E(1)_{n}(X)=M_{n}$ defined in 1.1 and by taking as Adams operations, the Adams operations just mentioned.

\vspace{30pt}
Remarks 1.

\vspace{12pt}

a) We know from ([B83]page 929) that  $\cal U$: $\pi_{*}(E(1)-Mod \mapsto $ $\cal B$ ${(l)}_{*}$ verifies:

\vspace{5pt}
$\cal U$ $(G)= E(1)_{*}E(1)\otimes_{\pi_{*}E(1)} G$

\vspace{5pt}

for all $\pi_{*}(E(1))$-module $G$. In particular taking $0$ component, $\cal U$ $_{0}(G)=E_{0}(1)E(1)\otimes_{Z_{(l)}} G $

\vspace{6pt}

 Therefore if $\phi: G \mapsto G$ is a map of $Z_{(l)}$-modules with an eigenvalue $\lambda$,

 \vspace{5pt}

 then $\cal U$ $ (\phi)$ has also eigenvalue $\lambda$ in $\cal U$ $(G)$.

\vspace{12pt}

b) The next proposition is also proven in [T89]

\vspace{16pt}

Proposition 1:\emph{Under the hypothesis that $\pi_{0}(L_{E(1)}(K(X_{\infty})))$ is $l$-reduced we get:}

\vspace{12pt}

\[\pi_{0}([L_{E(1)}K(X_{\infty})]^{l})\otimes Q \simeq [\pi_{0}(L_{E(1)}K(X_{\infty}))]^{l}\otimes Q\]

\vspace{12pt}

Proof: Let $G_{\nu}= \pi_{0}(Y\wedge M(Z/l^{\nu}))=\pi_{0}(T_{\nu})=\pi_{0}(Y/l^{\nu})$, where $Y=L_{E(1)}K(X_{\infty})$.There is an exact sequence  ([BK72] Chap 9)

\vspace{8pt}

(3)

\[0 \mapsto lim^{1}(G_{\nu}) \mapsto \pi_{0}(homlim T_{\nu}) \mapsto lim( G_{\nu}) \mapsto 0\]

\vspace{12pt}

 and (4)
  \[G_{\nu}= \pi_{0}(Y\wedge M(Z/l^{\nu})) \simeq (\pi_{0}(Y)\otimes Z/l^{\nu})\oplus Tor^{1}(Z/l^{\nu},Y)\]

  \vspace{8pt}

  Now $lim Tor^{1}(Z/l^{\nu},Y)=0$ since the limit is equal to

  \vspace{8pt}

\[{\Pi_{i=1}^{\infty}\{{ g_{i}=l^{i}-torsion-element\in \pi_{0}(L_{E(1)}K(X_{\infty}))}/ lg_{i+1}=g_{i}}\}\]

\vspace{12pt}

and each coordinate in this limit is $0$,for it belongs to the intersection of all $l^{l}(\pi_{0}L_{E(1)}K(X_{\infty}))$ which is $0$ because $\pi_{0}L_{E(1)}K(X_{\infty})$ is reduced.Then by (4), I get:

\vspace{12pt}

(5)
 \[lim( G_{\nu})= lim (\pi_{0}(Y)\otimes Z/l^{\nu})=[\pi_{0}(Y)]^{l}\].

Obviously (6): $lim^{1}(\pi_{0}(Y)\otimes Z/l^{\nu})=0$. On the other hand, $lim^{1}Tor^{1}(Z/l^{\nu},Y)$ has bounded $l$-torsion. Let me show why this is so:

\vspace{12pt}

 Let $M_{\nu}= Tor^{1}(Z/l^{\nu},Y)$.The map $M_{\nu+1} \mapsto M_{\nu}$ is the map which goes from the $l^{\nu+1}$-torsion elements of $Y$ to the $l^{\nu}$-torsion elements of $Y$ given by $ x \mapsto lx$.It is in general not surjective,so that it is difficult to prove that $lim^{1}M_{\nu}=0$. Anyway,(3) has simplified because of (5) and (6) to

\vspace{12pt}

 (7)
 \[0 \mapsto lim^{1}(Tor^{1}(Z/l^{\nu},\pi_{0}(Y)) \mapsto \pi_{0}(homlim Y\wedge M(Z/l^{\nu}))=\pi_{0}(Y^{l}) \mapsto lim( G_{\nu})=[\pi_{0}(Y)]^{l} \mapsto 0\]

 \vspace{8pt}

 where $Y=L_{E(1)}K(X_{\infty})$ and $lim( G_{\nu})=[\pi_{0}(Y)]^{l}$ is reduced since it is the projective limit of the reduced groups $ G_{\nu}$ See [MGM10]. $[\pi_{0}(Y)]^{l}$ is also a cotorsion group (see [F1]), since it is the epimorphic image of a cotorsion group in the exact sequence (7): $\pi_{0}(Y^{l})$ is equal to the cotorsion reduced group $Ext^{1}(Q/Z_{(l)},Y)$ since(see [B79])

 \vspace{12pt}

 \[Ext^{1}(Q/Z_{(l)},Y) \mapsto  \pi_{0}(Y^{l})\mapsto Hom(Q/Z_{(l)},\pi_{-1}(Y))\]

 \vspace{8pt}

 and $Hom(Q/Z_{(l)},\pi_{-1}(Y))=0$ since $\pi_{-1}(Y)$is $l$-reduced, see ([MGM10]) and therefore, $Ext^{1}(Q/Z_{(l)},Y) \simeq \pi_{0}(Y^{l})$ and $Ext^{1}(Q/Z_{(l)},Y)$ is a cotorsion group.The torsion group of the cotorsion reduced group $[\pi_{0}(Y)]^{l}$ noted $T([\pi_{0}(Y)]^{l})$ is in the terminology of [F2] an $l$-complete torsion group.Now $[\pi_{0}(Y)]^{l}$ is a reduced algebraically compact group since it is complete in the terminology of[R08] page 440.It is complete because it is the closure in the $l$-adic topology of the topological Hausdorff group $\pi_{0}(Y)]$ Being reduced and algebraically compact implies it is a direct summand of a direct product of cyclic $l$-groups by [F2] Corollary 38.2 page 161.Henceforth, $T([\pi_{0}(Y)]^{l})$ is contained in a direct product of cyclic $l$-groups, and so the torsion part of $\pi_{0}(Y)$, noted $T(\pi_{0}(Y))$ is contained in a direct sum of cyclic $l$-groups. Since $T(\pi_{0}(Y))$ is reduced, then it has bounded $l$-torsion, ie there exists $\nu_{0}$ such that $l^{\nu_{0}}T(\pi_{0}(Y)=0$ Then by definition of $lim^{1}$, $lim^{1}(Tor^{1}(Z/l^{\nu},\pi_{0}(Y))$ has bounded $l$-torsion, as wanted. Then,tensoring with $Q$ in the exact sequence (7) becomes,

 \vspace{12pt}

 (8)
 \[0 \mapsto \pi_{0}(Y^{l})\otimes Q \mapsto [\pi_{0}(Y)]^{l}\otimes Q \mapsto 0\]

 \vspace{8pt}

 and therefore Proposition 1 has been proved.

\vspace{16pt}

Remark 2: We conjecture that $\pi_{0}(Y^{l})$ and therefore also $[\pi_{0}(Y)]^{l}$ are without torsion, in which case by [R08] page 445 $[\pi_{0}(Y)]^{l}$, is a direct summand of copies of $Z_{l}$. Since tensored by $Q$, ie $[\pi_{0}(Y)]^{l}\otimes Q$, using (1) and the above Proposition 1, is a finite direct sum of copies of $Q_{l}$, $[\pi_{0}(Y)]^{l}$ has to be a finite direct sum of copies of $Z_{l}$.

\vspace{16pt}

Theorem 1: $\cal U$$([\pi_{0}(L_{E(1)}K(X_{\infty}))]^{l}\otimes Q)\simeq [E(1)_{0}(K(X_{\infty})]^{l}\otimes Q$

\vspace{16pt}

Corollary 2: $\cal U$$(\oplus_{i=0}^{d}H^{2i}_{et}(X_{\infty},Q_{(l)}(i)))\simeq [E(1)_{0}(K(X_{\infty})]^{l}\otimes Q$.

\vspace{12pt}

Proof:  $\cal U$$(\oplus_{i=0}^{d}H^{2i}_{et}(X_{\infty},Q_{(l)}(i)))\simeq $$\cal U$$([\pi_{0}(L_{E(1)}K(X_{\infty}))]^{l}\otimes Q)$.This fact follows from Proposition 1 and from (1). Hence, corollary 2 follows immediately from theorem 1.

\vspace{16pt}

Proof of theorem 1:

\vspace{12pt}

Since $\pi_{0}(L_{E(1)}K(X_{\infty}))$ is $l$-reduced the kernel of the $l$-completion
\vspace{2pt}
$\pi_{0}(L_{E(1)}K(X_{\infty})) \mapsto [\pi_{0}(L_{E( 1)}K(X_{\infty}))]^{l}$ is equal to $0$ and the completion map is injective. Also the $l$-adic topology in $\pi_{0}(L_{E(1)}K(X_{\infty}))$ is Hausdorff and this space is dense in its $l$-completed space. The functor $\cal U$ is exact and $0 \mapsto \cal U$$(\pi_{0}(L_{E(1)}K(X_{\infty}))\otimes Q)\simeq E(1)_{0}K(X_{\infty})\otimes Q \mapsto $$\cal U$$([\pi_{0}(L_{E(1)}K(X_{\infty}))]^{l})\otimes Q$ has dense image. On the other hand since $\pi_{0}(L_{E(1)}K(X_{\infty}))$ is $l$-reduced, $ \cal U$$(\pi_{0}(L_{E(1)}K(X_{\infty})))\otimes Q \simeq E(1)_{0}K(X_{\infty})\otimes Q$ is $l$-reduced (See [B83]) and the isomorphism holds as an isomorphism of $Q$-vector spaces.Then,

\vspace{8pt}

\[0 \mapsto E(1)_{0}K(X_{\infty})\otimes Q \mapsto [E(1)_{0}K(X_{\infty})]^{l}\otimes Q\]

\vspace{8pt}

with dense image. By uniqueness of the $l$-completed space, (9):$\cal U$$([\pi_{0}(L_{E(1)}K(X_{\infty}))]^{l})\otimes Q\simeq [E(1)_{0}(K(X_{\infty})]^{l}\otimes Q$, and the theorem follows.

\vspace{16pt}

Remark 3: If the conjecture stated in remark 2 holds, then by (9), remark 2, and remark 1 b), $[E(1)_{0}(K(X_{\infty})]^{l}$ is a finite direct sum of copies of $Z_{l}[[t]]$, a fact which was proved for nonsingular complete curves in [DM95].

\vspace{16pt}

Theorem 2: \emph{ The geometric frobenius $\Phi_{X}$ acts semisimply in $\oplus_{i=0}^{d} H^{2i}_{et}(X_{\infty},Q_{l}(i))$}.

\vspace{12pt}

Before we give a proof we need to state one definition and two remarks:

\vspace{12pt}

Definition 1: Let $E$ and $F$ be cohomology theories. A natural transformation from $E^{n}(-)$ to $F^{m}(-)$ is called a cohomology operation from $E^{n}(-)$ to $F^{m}(-)$. If it is compatible with the suspension isomorphisms then it is called a stable operation. (See [Bo95])

\vspace{12pt}

Remark 4: Let $E$ and $F$ be spectra. Then the set of stable cohomology operations from $E$ to $F$ can be identified with $F^{*}(E)$ [Bo95].Therefore the ring $E(1)^{*}E(1)$ may be identified with the stable operations of degree $0$, $\phi: E(1)_{*}(-) \mapsto E(1)_{*}(-)$ which in turn are induced by map of spectra $\phi: E(1) \mapsto E(1)$ (See [KJ84] page 57) Therefore given a base for $E(1)^{*}E(1)$, we obtain a base for the stable operations of degree $0$ on $E(1)_{*}(-)$.

\vspace{12pt}

Remark 5: From ([CCW05] page 13), we know that  $\widehat{E(1)}E(1)$ is isomorphic to $Z_{l}[[Y]]$, where $Y=\Psi^{r}-1$ for the Adams operation $\Psi^{r}$ with $r$ a primitive modulo $l^2$ and where $\widehat{E(1)}$ is the $l$-adic completion of $E(1)$. It is in particular from this fact  that Bousfield obtains in ([B83]page 908) an equivalence between the category $\cal B$ $(l)$ and the category $\cal B$ $(l)^{r}$. This isomorphism gives us a base for the ring $\widehat{E(1)}E(1)$, and by remark 4, it gives us a base for the stable degree $0$ operations $\phi: E(1)_{*}(-) \mapsto E(1)_{*}(-)$, and in particular for the $0$ component degree $0$ operations $\phi: E(1)_{0}(-) \mapsto E(1)_{0}(-)$

\vspace{12pt}

Proof of theorem 2:
\vspace{12pt}

By the Weil Conjectures the eigenvalues of $\Phi_{X}$ acting on $H^{2i}_{et}(X_{\infty},Q_{l}(i))$  are algebraic numbers all of whose complex conjugates have absolute value $q^{i}$.
\vspace{3pt}

Then $\cal U$$(\Phi_{X})$ acting on $\cal U$$(\oplus_{i=0}^{d}H^{2i}_{et}(X_{\infty},Q_{l}(i)))$ has eigenvalues whose complex conjugates have absolute value $q^{i}$,$ i\in {1,2,..d}$.

\vspace{12pt}

We will prove in a moment that this map can be identified with the Adams operation $\Psi^{q}$ on $[E(1)_{0}K(X_{\infty})]^{l}$ which is an object in $\cal B$. Now, this Adams operation, $\Psi^{q}$, is diagonalizable  on $[E(1)_{0}K(X_{\infty})]^{l}\otimes Q $ Then, by corollary 2, it is diagonalizable in $\cal U$$(\oplus_{i=0}^{d}H^{2i}_{et}(X_{\infty},Q_{l}))$ . This fact in turn implies that $\cal U$$(\Phi_{X})$ is also diagonalizable in $\cal U$$(\oplus_{i=0}^{d}H^{2i}_{et}(X_{\infty},Q_{l}))$, which then implies that $\Phi_{X}$ is diagonalizable in $\oplus_{i=0}^{d}H^{2i}_{et}(X_{\infty},Q_{l})$ as wanted.

\vspace{12pt}

 $\cal U$$(\Phi_{X})$ can be identified with the Adams operation $\Psi^{q}$ because by the isomorphism of remark 5, $\cal U$$(\Phi_{X})$ is an infinite combination of the elements,$(\Psi^{r}-1)^{s}$ $s=0,1,2,3...$. This last fact says that the eigenvalues of $\cal U$$(\Phi_{X})$ on $[E(1)_{0}K(X_{\infty})]^{l}\otimes Q $ are a combination of the of the eigenvalues of $(\Psi^{r}-1)^{s}$, $s=0,1,2,3...$ on the same space, which are real numbers, implying that the eigenvalues of $\cal U$$(\Phi_{X})$ are real numbers. Since on the other hand they are algebraic numbers whose complex conjugate have absolute value $q^{i}$, they must be equal to $q^{i}$ and hence it is the Adams operation $\Psi^{q}$, as stated above.

\vspace{16pt}

\emph{E-mail address}: valmont8ar@hotmail.com

\end{document}